\documentclass[11pt,A4]{amsart}
\usepackage{ifthen,latexsym,amssymb,amsmath}
\newtheorem{theorem}{Theorem}
\newtheorem{lemma}{Lemma}

\newtheorem{conjecture}{Conjecture}

\begin{document}

\title{On the height of cyclotomic polynomials}

\author{Bart\l{}omiej Bzd\c{e}ga}

\address{Str\'o\.zy\'nskiego 15A/20 \\ 60-688 Pozna\'n, Poland}

\email{exul@wp.pl}



\maketitle

\begin{abstract}
Let $A_n$ denote the height of cyclotomic polynomial $\Phi_n$, where $n$ is a product of $k$ distinct odd primes. We prove that $A_n \le \varepsilon_k\varphi(n)^{k^{-1}2^{k-1}-1}$ with $-\log\varepsilon_k\sim c2^k$, $c>0$. The same statement is true for the height $C_n$ of the inverse cyclotomic polynomial $\Psi_n$.

Additionally, we improve on a bound of Kaplan for the maximal height of divisors of $x^n-1$, denoted by $B_n$. We show that $B_n<\eta_kn^{(3^k-1)/(2k)-1}$, with $-\log \eta_k \sim c3^k$ and the same $c$.
\end{abstract}

\section{Introduction}

The polynomial
$$\Phi_n(x)=\sum_{0\le m\le \varphi(n)}a_n(m)x^m\quad=\prod_{k\le n,\;(k,n)=1}(x-\zeta_n^k)$$
where $\zeta_n=e^{2i\pi/n}$, is called the $n$-th cyclotomic polynomial. We are interested in estimating its coefficients, so we define
$$A_n=\max_m|a_n(m)| \quad \text{and} \quad S_n=\sum_{m=0}^{\varphi(n)}|a_n(m)|.$$
We define also
$$\Psi_n(x)=\frac{1}{\Phi_n(x)}=\sum_{m\ge0}c_n(m)x^m, \quad C_n=\max_m |c_n(m)|.$$
The polynomial $(1-x^n)\Psi_n(x)$ is called the $n$-th inverse cyclotomic polynomial (see \cite{Moree} for details). We remark that $c_n(m)$ is equal to the $m'-$th coefficient of the $n-$th inverse cyclotomic polynomial, where $0\le m' < n$ and $m'\equiv m \pmod n$.

We consider the numbers $n$ which are odd and square free only, since it is known that $A_{\ker(n)}=A_n=A_{2n}$, where $\ker(n)$ is the product of all distinct prime factors of $n$ (see \cite{Than} for details). The same fact is true for inverse cyclotomic polynomials.

The order of $\Phi_n$ is the number $\omega(n)$ of primes dividing $n$. For $\omega(n)\le4$ the following bounds are known:
\begin{equation} \label{1234}
A_p=1, \quad A_{pq}=1, \quad A_{pqr} \le \epsilon_3p, \quad A_{pqrs} \le \epsilon_4p^3q.
\end{equation}
The first of them is obvious. The second one is due to A. Migotti \cite{Mig}.

The third one with $\epsilon_3=1$ is due to A. S. Bang \cite{Bang}. It has been improved by some authors. Presently it is known that one can take $\epsilon_3=3/4$ (see  \cite{Bach,Beit,Bzd2}) and that one cannot replace $\epsilon_3$ by a constant smaller than $2/3$ (see \cite{GM}). It is strongly believed that the estimate holds with $\epsilon_3=2/3$ (J. Zhao and X. Zhang \cite{ZZ}, preprint). This conjecture is known as the Corrected Beiter Conjecture (see \cite{GM}).

The fourth inequality with $\epsilon_4=1$ was established by Bloom \cite{Bloom}. We use a simple argument from \cite{BPV} to show that the inequality is true with $\epsilon_4=\epsilon_3$.

For inverse cyclotomic polynomials we know the following bounds
$$C_p=1, \quad C_{pq}=1, \quad C_{pqr}\le p-1.$$
The first and the second of them are easy to obtain. The third was proved by P. Moree \cite{Moree} who in the same paper proved that $p-1$ cannot be replaced by a smaller number.

In the general case, we know the following result by  P. T. Bateman, C. Pomerance and R. C. Vaughan \cite{BPV} for standard cyclotomic polynomials.

\begin{equation} \label{old}
A_{p_1\ldots p_k} \le M_k \le n^{k^{-1}2^{k-1} - 1},
\end{equation}
where $M_k=\prod_{i=1}^{k-2}p_i^{2^{k-i-1}-1}$ (this notation we use troughout the paper). The same authors came up with the following conjecture (cf. \cite{BPV}, p. 175).

\begin{conjecture} \label{BPV con}
In (\ref{old}) one can replace $n$ by $\varphi(n)$.
\end{conjecture}

We prove this conjecture and moreover, we improve it by multiplying the right hand side by a constant depending on $k$ only and decreasing quickly when $k$ grows. We prove also a similar result for the inverse cyclotomic polynomials and give the bound for the maximal magnitude of the coefficient of any divisor of $x^n-1$, improving on an earlier result of N. Kaplan \cite{Kap}. The idea of estimating the maximal magnitude of coefficient of any divisor of $x^n-1$ comes from C. Pomerance and N. C. Ryan \cite{PR}.

By $\epsilon_k$ we denote the smallest positive real number for which the inequality $A_{p_1\ldots p_k} \le \epsilon_kM_k$ holds with any distinct primes $p_1,\ldots,p_k$. In the same way we define $\epsilon_k^{inv}$ for the inverse cyclotomic polynomial and $E_k$. Let
\begin{equation}\label{def}
d=\max_{p,q,r} \frac{S_{pqr}}{p^2qr}, \; \rho=\prod_{i=0}^{\infty}\left({\frac{2i+5}{2i+6}}\right)^{1/2^i}, \; C=\left(\frac34\epsilon_3^{3/2}d\rho^{1/4}\right)^{1/32}.
\end{equation}
Note that $C<1$. Our main results are the four following theorems.

\begin{theorem} \label{main}
We have $\log \epsilon_k \le 2^k\log(C+o(1))$
\end{theorem}

\begin{theorem} \label{inv}
We have $\log\epsilon_k^{inv} \le 2^k\log(C+o(1))$
\end{theorem}

\begin{theorem} \label{Kaplan+}
If $B_n = \eta_kn^{(3^k-1)/(2k)-1}$, then $\log \eta_k < 3^k\log(C+o(1))$ for every $n$ free of squares.
\end{theorem}

\begin{theorem} \label{phi}
Conjecture \ref{BPV con} holds true, that is we have $M_k\le\varphi(n)^{k^{-1}2^{k-1}-1}$ with $n=p_1\ldots p_k$.
\end{theorem}

In the proof of Theorem \ref{main} we also establish the following bounds
\begin{equation} \label{small k}
A_{pqrs}\le \frac{3}{4}p^3q, \quad A_{pqrst}\le \frac{135}{512}p^7q^3r, \quad A_{pqrstu}\le \frac{18225}{262144}p^{15}q^7r^3s,
\end{equation}
where we assumed $\epsilon_3=3/4$. For $\epsilon_3=2/3$ we establish constants $\frac{2}{3}$, $\frac{2}{9}$, $\frac{32}{729}$, respectively.

Also for the inverse cyclotomic polynomial
\begin{equation} \label{small k inv}
C_{pqrs} \le \frac34p^3q, \quad C_{pqrst} \le \frac{9}{16}p^7q^3r, \quad C_{pqrstu} \le \frac{10935}{131072}p^{15}q^7r^3s
\end{equation}
for $\epsilon_3=3/4$. If $\epsilon_3=2/3$, then we obtain constants $\frac{2}{3}$, $\frac{4}{9}$, $\frac{8}{81}$, respectively

Let us remark that Theorem \ref{main}, but with larger constant, can be obtained by the original method of P. T. Bateman, C. Pomerance and R. C. Vaughan. Our method is a bit different. It is based on a different recursive formula given in Lemma \ref{tool}. We use also some basic combinatorics.

\section{Preliminaries}

Our primary tool is the following lemma.

\begin{lemma} \label{tool} Let $p_1,\ldots,p_k$ be distinct primes. Then
\begin{equation} \label{l}
\Phi_{p_1\ldots p_k}(x)=f(x) \cdot \prod_{j=1}^{k-2}P_j(x),
\end{equation}
where
\begin{equation} \label{f}
f(x)=(1-x^{p_1\ldots p_k}) \cdot \frac{\prod_{i=2}^k(1-x^{p_2\ldots p_k/p_i})}{\prod_{i=1}^k(1-x^{p_1\ldots p_k/p_i})}
\end{equation}
and $P_j=\prod_{i=j+2}^k\Phi_{p_1\ldots p_j}(x^{p_{j+2}\ldots p_k/p_i})$.
\end{lemma}
As $\deg(\Phi_n)=\varphi(n)<n$, we may replace $f$ by $f^* \equiv f \pmod {x^{p_1\ldots p_k}}$ in (\ref{l}), where $\deg(f^*)<p_1\ldots p_k$. Then we have congruence modulo $x^{p_1\ldots p_k}$ in (\ref{l}) instead of equality, which does not matter for our purposes. In addition in the next section we prove the following lemma.
\begin{lemma} \label{f*}
We have $H(f^*) \le b_{k-2}=\binom{k-2}{\lfloor (k-2)/2 \rfloor}$.
\end{lemma}

Lemmas \ref{tool} and \ref{f*} allow us to give the following recursive bound on $\epsilon_k$.

\begin{lemma} \label{epsilon2}
We have $\epsilon_k \le E_k = \frac{b_{k-2}d^{k-4}}{2^{k-3}}\prod_{j=1}^{k-2}\epsilon_j^{k-j-1}$.
\end{lemma}

To start the induction we need also the following estimates.

\begin{lemma} \label{epsilon1}
We have $\epsilon_4 \le \epsilon_3$.
\end{lemma}

\begin{proof}
It is known that $S_1=2$ and $S_{pq}\le pq/2$ (see \cite{Bloom} for a proof of the second equality). By Lemma 4 on pages 182--183 in \cite{BPV},
$$A_{pqrs} \le A_{pqr}S_{pq}S_pS_1 \le \epsilon_3 \cdot p^3q,$$
so the estimate holds.
\end{proof}

\begin{lemma} \label{d}
For $d$ defined in (\ref{def}) we have $d \le \epsilon_3(2-\epsilon_3)/2$.
\end{lemma}

\begin{proof}
Bloom \cite{Bloom} proved that
$$|a_{pqr}(m)|=|a_{pqr}(\varphi(pqr)-m)| \le 2(\lfloor m/qr \rfloor + 1).$$
Thus
\begin{eqnarray*}
S_{pqr} & \le & 2\sum_{k=0}^{\varphi(pqr)/2}\min\{\epsilon_3p,2(\lfloor m/qr \rfloor + 1)\} \\
& \le & \epsilon_3p(\varphi(pqr)+2 - 2\lfloor \epsilon_3p/2 \rfloor qr) + 2qr\sum_{a=0}^{\lfloor \epsilon_3p/2 \rfloor-1}(2a+2) \\
& = & \epsilon_3p(p-1)(q-1)(r-1) + 2\epsilon_3p - 2\lfloor \epsilon_3p/2 \rfloor \epsilon_3pqr \\
&& + 2\lfloor \epsilon_3p/2 \rfloor(2\lfloor \epsilon_3p/2 \rfloor+1)qr \\
& < & \epsilon_3(2-\epsilon_3)p^2qr/2,
\end{eqnarray*}
which completes the proof.
\end{proof}

\bigskip\section{Proof of Lemma \ref{tool}, \ref{f*} and \ref{epsilon2}}

\begin{proof}[Proof of Lemma \ref{tool}]
We prove this lemma by induction on $k$. By (see \cite{Bloom}) it holds for $k<5$. Let us define
$$\widetilde{f}(x)=(1-x^{p_2\ldots p_k}) \cdot \frac{\prod_{i=3}^k(1-x^{p_3\ldots p_k/p_i})}{\prod_{i=2}^k(1-x^{p_2\ldots p_k/p_i})}$$
and $\widetilde{P}_j(x)=\prod_{i=j+2}^k\Phi_{p_2\ldots p_j}(x^{p_{j+2}\ldots p_k/p_i})$. By the inductive assumption,
\begin{equation} \label{ind}
\Phi_{p_2\ldots p_k}=\widetilde{f}(x) \cdot \prod_{j=2}^{k-2}\widetilde{P}_j(x).
\end{equation}
It is known that $\Phi_{np}(x) = \Phi_n(x^p) / \Phi_n(x)$ for a prime $p$ not dividing $n$ (see \cite{Than}). Then also
$$\Phi_{p_1\ldots p_k}(x)=\frac{\Phi_{p_2\ldots p_k}(x^{p_1})}{\Phi_{p_2\ldots p_k}(x)} \quad \text{and} \quad P_j(x)=\frac{\widetilde{P}_j(x^{p_1})}{\widetilde{P}_j(x)}.$$
By this and (\ref{ind})
$$\Phi_{p_1 \ldots p_k}(x) = \frac{\widetilde{f}_k(x^{p_1}) \cdot \prod_{j=2}^{k-2}\widetilde{P}_j(x^{p_1})}{\widetilde{f}_k(x) \cdot \prod_{j=2}^{k-2}\widetilde{P}_j(x)}
= \frac{\widetilde{f}(x^p)}{\widetilde{f}(x)P_1(x)} \cdot \prod_{j=1}^{k-2}P_j(x).$$
Finally,
$$\frac{\widetilde{f}(x^{p_1})}{\widetilde{f}(x)}=P_1(x)(1-x^{p_1\ldots p_k})\cdot \frac{\prod_{i=2}^k(1-x^{p_2\ldots p_k/p_i})}{\prod_{i=1}^k(1-x^{p_1\ldots p_k/p_i})}=P_1(x)f(x),$$
which completes the proof.
\end{proof}

\begin{proof}[Proof of Lemma \ref{f*}]
Let $n=p_1\ldots p_k$ and $f^*(x)=\sum_{m=0}^{n-1}d_mx^m$. By (\ref{f}) we have
\begin{equation} \label{dm1}
f^*(x) \equiv \prod_{i=2}^k(1-x^{p_2\ldots p_k/p_i}) \sum_{\alpha_1,\ldots,\alpha_k\ge0}x^{\alpha_1n/p_1+\ldots+\alpha_kn/p_k} \pmod{x^n}.
\end{equation}
Let
$$\Lambda=\{\lambda=(\lambda_2,\ldots,\lambda_k):\lambda_i\in\{0,1\}\text{ for }i=2,\ldots,k\}, \quad s(\lambda)=(-1)^{\lambda_2+\ldots+\lambda_k}.$$
By (\ref{dm1})
\begin{equation} \label{dm2}
d_m=\sum_{\lambda\in\Lambda}s(\lambda)\chi(m-\langle\lambda,v/p_1\rangle),
\end{equation}
where $\langle\cdot,\cdot\rangle$ is the scalar product in $\mathbb{R}^{k-1}$, $v = (n/p_2,\ldots,n/p_k)$ and
$$\chi(m) = \left\{\begin{array}{ll}
1 & \text{if } m \text{ is of the form } \alpha_1n/p_1+\ldots+\alpha_kn/p_k, \\
0 & \text{otherwise}.
\end{array}\right.$$
We define a number $\beta(\lambda)$ and a vector $\alpha(\lambda)=(\alpha_2(\lambda_2),\ldots,a_k(\lambda_k))$ by the congruence
\begin{equation} \label{m-<>}
m-\langle\lambda,v/p_1\rangle \equiv \beta(\lambda) n/p_1 + \langle\alpha(\lambda),v\rangle \pmod n.
\end{equation}
The numbers $\alpha_i(0)$ and $\alpha_i(1)$ depend only on the residue class of $m$ modulo $p_i$, so (\ref{m-<>}) holds for every $\lambda\in\Lambda$. We have the following equivalences
\begin{eqnarray*}
& & \chi(m-\langle\lambda,v/p_1\rangle)=1 \\
& \iff & \langle\lambda,v/p_1\rangle + \langle\alpha(\lambda),v\rangle \le m \\
& \iff & \langle\lambda,v/p_1\rangle + \langle\alpha(\lambda)-\alpha(\theta_{k-1}),v\rangle \le m - \langle\alpha(\theta_{k-1}),v\rangle, \\
\end{eqnarray*}
where $\theta_{k-1}=(0,\ldots,0)$. We have
$$\langle\alpha(\lambda)-\alpha(\theta_{k-1}),v\rangle = \sum_{i=2}^k(\alpha_i(\lambda_i)-\alpha_i(0))v_i=\sum_{i=2}^k(\alpha_i(1)-\alpha_i(0))v_i\lambda_i = \langle\lambda,w\rangle,$$
where $w=((\alpha_i(1)-\alpha_i(0))v_i)_{i=2}^k$. Therefore
$$\chi(m-\langle\lambda,v/p_1\rangle)=1 \iff \langle\lambda,u\rangle \le D,$$
where $u=v/p_1 + w$ and $D=m-\langle\alpha(\theta_{k-1}),v\rangle$. By (\ref{dm2})
\begin{equation} \label{dm3}
d_m= \sum_{\lambda\in\Lambda,\;\;\langle\lambda,u\rangle \le D}s(\lambda).
\end{equation}

Without loss of generality we may assume that $0 \le u_k \le u_2,\ldots,u_{k-1}$.

There is a natural bijection between $\Lambda$ and the family of subsets of $\{2,3,\ldots,k\}$, defined by
$$S_{\lambda}=\{i\in\{2,\ldots,k\}:\lambda_i=1\} \quad\text{for } \lambda \in \Lambda.$$
We say that $\lambda=(\lambda_2,\ldots,\lambda_{k-1},0)$ is maximal if $\langle\lambda,u\rangle \le D$ and for every $\lambda'=(\lambda'_2,\ldots,\lambda'_{k-1},0)$ such that $S_{\lambda}\subset S_{\lambda'}$ we have $\langle\lambda',u\rangle > D$. Note that for
$$\lambda^0=(\lambda_2,\ldots,\lambda_{k-1},0) \quad \text{and} \quad \lambda^1=(\lambda_2,\ldots,\lambda_{k-1},1)$$
the following statements are true.
\begin{itemize}
\item If $\lambda^0$ is not maximal and $\langle\lambda^0,u\rangle \le D$ then $\langle\lambda^1,u\rangle \le D.$
\item If $\langle\lambda^1,u\rangle \le D$ then $\langle\lambda^0,u\rangle \le D$.
\item $s(\lambda^0) + s(\lambda^1) = 0$.
\end{itemize}
By this observation and (\ref{dm3}) we conclude that
\begin{equation} \label{dm4}
|d_m| \le \#\{\lambda\in\Lambda \; : \; \lambda \text{ is maximal}\}.
\end{equation}
Let $\lambda^1,\ldots,\lambda^t \in \Lambda$ be maximal. By the definition of maximal $\lambda$, we have $S_{\lambda^i}\subset\{2,\ldots,k-1\}$ and $S_{\lambda^i}\not\subset S_{\lambda^j}$ for every $i\neq j$.
\begin{theorem}[E. Sperner, 1928]\label{Sperner}
Let $A_1,\ldots,A_t \subset A$, where $\#A\le\infty$. If $A_i\not\subset A_j$ for every $i\neq j$, then $t\le\binom{\#A}{\lfloor \#A/2 \rfloor}$. \qed
\end{theorem}
For the proof see \cite{Sper}.

By Theorem \ref{Sperner} and (\ref{dm4}), $|d_m|\le t\le\binom{k-2}{\lfloor(k-2)/2\rfloor}$.
\end{proof}

\begin{proof}[Proof of Lemma \ref{epsilon2}]
For a formal power series $f(x)=\sum_{m\ge0}a_mx^m \in \mathbb{Z}[[x]]$ we define $H,S\in[0,\infty]$
$$H(f)=\max_{m\ge0}|a_m|, \quad S(f)=\sum_{m\ge0}|a_m|.$$
We call $H(f)$ the height of $f$. Note that
\begin{equation} \label{H}
H\left(f(x)\prod_{i=1}^kQ_i(x)\right)\le H(f)\prod_{i=1}^kS(Q_i),
\end{equation}
\begin{equation} \label{S}
S\left(\prod_{i=1}^kQ_i(x)\right)\le \prod_{i=1}^kS(Q_i)
\end{equation}
for polynomials $Q_1,Q_2,\ldots,Q_{k} \in \mathbb{Z}[x]$ and a formal power series $f$. By (\ref{S}) we have for $j<k$
$$S_{p_1\ldots p_j} \le (\deg(\Phi_{p_1\ldots p_j})+1)A_{p_1\ldots p_j} \le \epsilon_j \cdot p_j \cdot p_1^{2^{j-2}}p_2^{2^{j-3}}\ldots p_{j-2}^2p_{j-1},$$
as $\deg(\Phi_n)=\varphi(n)<n$ for $n>1$. Then again by (\ref{S})
\begin{equation} \label{SP}
S(P_j) \le \epsilon_j^{k-j-1}\left(p_j \cdot p_1^{2^{j-2}}p_2^{2^{j-3}}\ldots p_{j-2}^2p_{j-1}\right)^{k-j-1},
\end{equation}
where $P_j$ is defined in Lemma \ref{tool}. Additionally,
\begin{equation} \label{SP23}
S_{p_1p_2}<p_1p_2/2, \quad S_{p_1p_2p_3} \le d \cdot p_1^2p_2p_3.
\end{equation}
Applying (\ref{H}), (\ref{SP}), (\ref{SP23}) and Lemma \ref{f*} to Lemma \ref{tool} we receive
\begin{eqnarray*}
A_{p_1\ldots p_k} & \le & \frac{b_{k-2}d^{k-4}}{2^{k-3}} \cdot \prod_{j=1}^{k-2}\epsilon_j^{k-j-1} \cdot \prod_{j=1}^{k-2}\left(p_j \cdot p_1^{2^{j-2}}p_2^{2^{j-3}}\ldots p_{j-2}^2p_{j-1}\right)^{k-j-1} \\
& = & E_kM_k,
\end{eqnarray*}
which completes the proof.
\end{proof}

\section{Proof of Theorem \ref{main}, \ref{inv}, \ref{Kaplan+} and \ref{phi}}

\begin{proof}[Proof of Theorem \ref{main}]
Consider a sequence $(e)$ given by the following conditions:
$$e_1=e_2=1, \quad e_3=e_4=\epsilon_3,$$
$$e_{k}=\frac{b_{k-2}d^{k-4}}{2^{k-3}}\prod_{j=1}^{k-2}e_j^{k-j-1} \quad \text{for } k\ge5.$$
By Lemmas \ref{epsilon2} and \ref{epsilon1} we have $\epsilon_k \le e_k$. We can easily compute that
\begin{equation} \label{small k proof}
e_5=\frac34\epsilon_3d, \quad e_6=\frac{9}{16}\epsilon_3^3d^2, \quad\ldots
\end{equation}
For $k\ge7$
$$\frac{e_k/e_{k-1}}{e_{k-1}/e_{k-2}}=e_{k-2}\cdot\frac{b_{k-2}b_{k-4}}{b_{k-3}^2},$$
then
$$e_k=e_{k-1}^2\cdot\frac{b_{k-2}b_{k-4}}{b_{k-3}^2},$$
therefore
$$e_k = e_6^{2^{k-6}}\cdot\prod_{i=7}^{k}\left(\frac{b_{i-2}b_{i-4}}{b_{i-3}^2}\right)^{2^{k-i}}.$$
Note that
$$\frac{b_{i-2}b_{i-4}}{b_{i-3}^2} = \left\{\begin{array}{ll}
\frac{i-2}{i-1}, & \text{for odd } i\\
\frac{i-2}{i-3}, & \text{for even } i.
\end{array}\right.$$
Then
\begin{eqnarray*}
e_k &= & e_6^{2^{k-6}}\quad\cdot\left(\frac56\right)^{2^{k-7}}\cdot\left(\frac65\right)^{2^{k-8}}\cdot\left(\frac78\right)^{2^{k-9}}\cdot\left(\frac87\right)^{2^{k-10}}\cdot\ldots \\
& = & \left(\frac{9}{16}\epsilon_3^3d^2\right)^{2^{k-6}}\cdot(1+o(1))\prod_{i=4}^{\lfloor k/2 \rfloor}\left(\frac{2i-3}{2i-2}\right)^{2^{k-2i}} \\
& = & \left(\frac34\epsilon_3^{3/2}d\rho^{1/4}+o(1)\right)^{2^{k-5}},
\end{eqnarray*}
which completes the proof of the Theorem \ref{main}.
\end{proof}

Note that (\ref{small k proof}) implies the bounds from (\ref{small k}).

\begin{proof}[Proof of Theorem \ref{inv}]
By the well known formula $\Psi_{np}(x)=\Psi_n(x^p)\Phi_n(x)$ we have
$$c_{np}(m)=\prod_{i=1}^{\lfloor m/p \rfloor}c_n(k)a_n(m-kp).$$
We note that $a_n(t)=0$ for $t\not\in\{0,\ldots,\varphi(n)\}$, and therefore
$$C_{p_1\ldots p_k} \le \left(\left\lfloor\frac{\varphi(p_1\ldots p_{k-1})}{p_k}\right\rfloor+1\right)A_{p_1\ldots p_{k-1}}C_{p_1\ldots p_{k-1}} \le p_1\ldots p_{k-2} \cdot A_nC_n$$
for $k\ge2$. Thus
$$C_{p_1\ldots p_k} \le C_{p_1p_2} \prod_{j=2}^{k-1}(p_1\ldots p_{j-1} \cdot A_{p_1\ldots p_j}) \le \epsilon_2\ldots \epsilon_{k-1}M_k.$$
Therefore
$$\epsilon_k^{inv} \le \epsilon_2\ldots \epsilon_{k-1} \le e_2\ldots e_{k-1}=\frac{b_{k-2}}{b_{k-3}}e_k$$
for $k\ge6$. It completes the proof.
\end{proof}

We can also prove that
$$\epsilon_4^{inv} \le \epsilon_3, \quad \epsilon_5^{inv} \le \epsilon_3^2, \quad \epsilon_6^{inv} \le \frac34\epsilon_3^3d$$
to justify (\ref{small k inv}).

\begin{proof}[Proof of Theorem \ref{Kaplan+}]
We recall that every divisor of $x^n-1$ is of the form $\prod_{d \in D}\Phi_d(x)$, where $D$ is a subset of the set of divisors of $n$. By (\ref{H}) and Theorem \ref{main}
\begin{eqnarray*}
B_n & \le & A_n\prod_{d\mid n,\;\;d<n}S_d \le \frac{2}{n}\prod_{d\mid n}dA_d \\
& \le & \frac{2}{n}\left(\prod_{d\mid n}d\right) \left(\prod_{d\mid n}\epsilon_{\omega(d)}\right) \left(\prod_{d\mid n}M_k(d)\right),
\end{eqnarray*}
where $M_k(d) = \prod_{i=1}^{\kappa-2}p_{\delta_i}^{2^{\kappa-i-1}-1}$ for $d=p_{\delta_1}\ldots p_{\delta_{\kappa}}$, $p_{\delta_1} < \ldots < p_{\delta_{\kappa}}$. We have
\begin{eqnarray*}
\frac{1}{n}\prod_{d\mid n}d & = & n^{2^{k-1}-1}, \\
\prod_{d\mid n}M_k(d) & \le & \prod_{\omega=1}^k \left(\left(\left(\sqrt[k]{n}\right)^{\omega}\right)^{2^{\omega-1}/\omega-1}\right)^{\binom{k}{\omega}} = n^{(3^k-1)/(2k)-2^{k-1}}.
\end{eqnarray*}
and by Theorem \ref{main}
\begin{eqnarray*}
\log\left(2\prod_{d\mid n}\epsilon_{\omega(d)}\right) & \le & \log2 + \sum_{d\mid n}2^{\omega(d)}\log(C+\xi_{\omega(d)}')\\
& \sim & 3^k\log C + \sum_{\omega=0}^k \binom{k}{\omega}2^{\omega}\xi_{\omega},
\end{eqnarray*}
where $\xi_{\omega}',\xi_{\omega}\to0$ with $\omega \to \infty$. It remains to prove that the sum equals $o(3^k)$. Indeed,
\begin{eqnarray*}
\sum_{\omega=0}^k \binom{k}{\omega}2^{\omega}\xi_{\omega} & \le & \xi_0\sum_{\omega=0}^{\lfloor \log k \rfloor}\binom{k}{\omega}2^{\omega} + \xi_{\lceil \log k \rceil}\sum_{\omega=0}^k\binom{k}{\omega}2^\omega \\
& = & O(2^{\log k}e^{\log^2 k}\log k) + o(3^k) = o(3^k),
\end{eqnarray*}
and the proof is done.
\end{proof}

In case $n=p_1\ldots p_k$ and $p_i\not\gg p_{i-1}$ for $i=2,\ldots, k$ Theorem \ref{Kaplan+} improves the result of N. Kaplan \cite{Kap} showing that
$$B_n < \prod_{j=1}^{k-1}p_j^{4\cdot3^{k-2}-1} \le n^{(4\cdot3^{k-2}-1)(k-1)/k}.$$

\begin{proof}[Proof of Theorem \ref{phi}]
We have $M_1=M_2=1$, so theorem holds for $k=1,2$. We prove it by induction on $k$. We assume that $p_1<\ldots<p_k$. Then for $k\ge3$
\begin{eqnarray*}
M_k & \le & p_1^{2^{k-2}-1} \cdot \varphi(p_2\ldots p_k)^{2^{k-2}/(k-1)-1} \\
& = & \left(\frac{p_1}{p_1-1}\right)^{\frac{2^{k-1}}{k}-1} \cdot \left(\frac{p_1^{k-1}}{\varphi(p_2\ldots p_k)}\right)^{\frac{2^{k-2}}{k-1}-\frac{2^{k-1}}{k(k-1)}} \cdot (\varphi(p_1\ldots p_k))^{\frac{2^{k-1}}{k}-1} \\
& \le & \left(\frac{p_1+1}{p_1}\right)^{\frac{2^{k-1}}{k}-1}\cdot\left(\frac{p_1}{p_1+1}\right)^{2^{k-2}-\frac{2^{k-1}}{k}}\cdot(\varphi(p_1\ldots p_k))^{\frac{2^{k-1}}{k}-1} \\
& \le & (\varphi(p_1\ldots p_k))^{2^{k-1}/k-1},
\end{eqnarray*}
which completes the proof of Theorem \ref{phi}.
\end{proof}

\section{Concluding remarks}

We analyze the value of the constant $C$. It is proved that $\epsilon_3\in[2/3,3/4]$, however we do not know the exact value of $\epsilon_3$. Similarly, we can only estimate the value of $d$. In the case $\epsilon_3=3/4$ and if we have the equality in Lemmas \ref{epsilon1} and \ref{d}, then $C \approx 0.953$. If the Corrected Beiter Conjecture holds, then $C \approx 0.946$.

Let us remark, that there exist a constant $\epsilon>0$ such that for $C<\epsilon$ the bound from Theorem \ref{main} is false. Indeed, if $p_j$ is the $j-$th odd prime number for $j\ge1$, then
$$1 \le A_{p_1\ldots p_k} \le (C+o(1))^{2^k}M_k$$
and therefore
$$C+o(1) \ge M_k^{-2^k} = \prod_{j=1}^{\infty}p_j^{-2^{3-j}}+o(1).$$
Using the prime number theorem we easily obtain that the product is convergent to a positive constant.

Recall the following conjecture of P. T. Bateman, C. Pomerance and R. C. Vaughan \cite{BPV}.
\begin{conjecture}
For every $k$ there exist a constant $\epsilon_k'$ such that
$$A_n \ge \epsilon_k'n^{2^{k-1}/k-1}$$
for infinitely many cyclotomic polynomials $\Phi_n$ of order $k$.
\end{conjecture}
If the conjecture is true, one of the most interesting questions is whether the maximal $\epsilon_k'$ is of the form $(C'+o(1))^{2^k}$ for some constant $0<C'<1$.

\section*{Acknowledgments}
The author would like to thank Pieter Moree for his suggestions how to make the paper more interesting and making some corrections. The author would also like to thank Wojciech Gajda for his remarks on the paper.

\end{document}